\documentclass[11pt]{amsart}
\usepackage{amsmath,amsthm,amsfonts,amssymb,amscd,hyperref}
\usepackage{fancyhdr}

\newtheorem{theorem}{THEOREM}[section]
\newtheorem{corollary}[theorem]{Corollary}
\newtheorem{remark}[theorem]{Remark}
\newtheorem{example}[theorem]{Example}
\newtheorem{lemma}[theorem]{Lemma}
\newtheorem{proposition}[theorem]{Proposition}
\newcommand{\nc}{\newcommand}
\newcommand{\s}{\sigma}
\newcommand{\g}{\mathcal{G}}
\newcommand{\al}{\alpha}
\newcommand{\ur}{\mathfrak{u}}
\nc{\scs}[1]{\scriptstyle{#1}}
\nc{\bin}[2]{  \left (\!\! \begin{array}{c} \scs{#1}\\
    \scs{#2} \end{array}\!\! \right )}  
\nc{\pr}{\partial}
\newcommand{\bt}{\begin{theorem}}
\newcommand{\et}{\end{theorem}}
\newcommand{\bco}{\begin{corollary}}
\newcommand{\eco}{\end{corollary}}
\newcommand{\bd}{\begin{definition}}
\newcommand{\ed}{\end{definition}}
\newcommand{\bex}{\begin{example}}
\newcommand{\eex}{\end{example}}
\newcommand{\bl}{\begin{lemma}}
\newcommand{\el}{\end{lemma}}
\newcommand{\bprop}{\begin{proposition}}
\newcommand{\eprop}{\end{proposition}}
\newcommand{\br}{\begin{remark}}
\newcommand{\er}{\end{remark}}
\newcommand{\bpf}{\begin{proof}}
\newcommand{\epf}{\end{proof}}
\title{Linear Differential Equations and Hurwitz Series}

\author{William F. Keigher\ \ }
\address{W.F. Keigher,  Department of Mathematics and Computer Science, Rutgers University, Newark, NJ 07102.}
\email{keigher@rutgers.edu}
\thanks{}

\author{\ \ V. Ravi Srinivasan}
\address{V. R. Srinivasan, Department of Mathematics and Computer Science, Rutgers University, Newark, NJ 07102.}
\email{ravisri@rutgers.edu}
\thanks{}

\dedicatory{This paper is dedicated to Professor Michael Singer on his sixtieth birthday}

\begin{document}
\maketitle

\begin{abstract}
In this article, we study the set of all solutions of linear differential
equations  using Hurwitz series.  We first obtain explicit
recursive expressions for solutions of such equations and study the group of
differential automorphisms of the set of all solutions. Moreover, we give explicit formulas that compute
the group of differential automorphisms. We require neither that the underlying field be algebraically closed nor that
the characteristic of the field be zero.
\end{abstract}

\section{Conventions and Basics}

Throughout, all rings are commutative with identity, and all
differential rings are ordinary (i.e., possess a single derivation,
which is often suppressed from the notation). Also, \textbf{N} will
denote the natural numbers $\{0,1,2,\ldots\}$ and \textbf{Q} the
field of rational numbers. Unless otherwise noted, $k$ will denote a
field.  If $V$ is a vector
space over $k$ and $X \subset V$, then ${\rm span}_kX$ will denote
the $k$-subspace of $V$ spanned by $X$. Let $R$ be a differential ring and let $y_1,y_2,\cdots,y_n\in R$. We denote the Wronskian of $y_1,y_2,\cdots,y_n$ by $w(y_1,y_2,\cdots,y_n)$. The set of all $n\times n$ matrices and $n\times n$ invertible matrices over
a field $k$ will be denoted by $M(n,k)$ and $GL(n,k)$ respectively. For $A\in M(n,k)$, we denote the centralizer of $A$ in $GL(n,k)$ by
$C_k(A):=\{T\in GL(n,k)|AT=TA\}$. Finally, for any $m,n \in
\mathbf{N}, \delta_n^m$ will denote the Kronecker delta, i.e.,
$\delta_n^m = 1$ if $m=n$ and $\delta_n^m = 0$ if $m \neq n$.

From~\cite{onthri} we recall that for any commutative ring $R$ with
identity, the \textsl{ring of Hurwitz series} over $R$, denoted by
$HR$, is defined as follows. The elements of $HR$ are sequences
$(a_n) = (a_0,a_1,a_2,\ldots)$, where $a_n \in R$ for each $n \in
\mathbf{N}$. Let $(a_n), (b_n) \in HR$.  Addition in $HR$ is defined
termwise, i.e.,
$$(a_n) + (b_n) = (c_n), \quad \text{where} \quad
c_n = a_n + b_n$$
for all $n \in \mathbf{N}$.  The (Hurwitz)
product of $(a_n)$ and $(b_n)$ is given by
$$(a_n)\cdot(b_n) = (c_n),  \quad
\text{where} \quad c_n = \sum_{j=0}^{n}{n \choose j}a_{j}b_{n-j}$$
for all $n \in \mathbf{N}.$ We
recall from~\cite{onthri} that if $\textbf{Q} \subseteq R$, then $HR
\cong R[[t]]$ via the mapping $(a_n) \mapsto
\sum_{n=0}^{\infty}\frac{a_n}{n!}t^n$.

Moreover, $HR$ is a differential ring with derivation
$\partial_{R}:HR \to HR$ given by
$$\partial_{R}((a_{0},a_{1},a_{2},\ldots)) =
(a_{1},a_{2},a_{3},\ldots).$$  We will often write $\partial$ in
place of $\partial_{R}$.  We have, as in~\cite{hufo}, for any $j \in
\mathbf{N}$, the additive mapping $\pi_j:HR \to R$ defined by
$\pi_j((a_n)) = a_j$.

In~\cite{onthri} it was shown that $H$ is a functor from
\textbf{Comm} (the category of commutative rings with identity) to
\textbf{Diff} (the category of ordinary differential rings) which is
the right adjoint to the functor $U:\mathbf{Diff} \to \mathbf{Comm}$
that ``forgets'' the derivation $d$ of a differential ring $(R,d)$.
This can be expressed as follows.

\bprop \label{rightadj}
    For any differential ring $(R,d)$ and any ring $S$, there is a
natural bijection between the sets of morphisms
\[\mathbf {Comm}(R,S) \cong
\mathbf {Diff}((R,d),(HS,\partial_{S})).\]  In particular, for any
ring homomorphism $f:R \to S$, there is a unique
differential ring homomorphism
$$\tilde{f}:(R,d) \to (HS,\partial_{S}) \quad
\text{given  by} \quad \tilde{f}(r) = (f(r), f(d(r)), f(d^{2}(r)),
\ldots).$$
\eprop

\section{Linear Homogenous Differential Operators}

Throughout this section, let $k$ be a field of any characteristic  and let $Hk$ be the differential ring of Hurwitz series over
$k$. Let $h_0,\ldots,h_{n-1} \in Hk$ and consider the monic linear
homogeneous differential operator
$$L: Hk \to Hk$$
defined for any $h \in Hk$ by
$$L(h) = \partial^{n}(h)+\sum_{i=0}^{n-1} h_{i}\partial^{i}(h).$$
We are interested in solutions to $L(h) = 0$ in $Hk$.  To this end,
let $V= \{ h \in Hk \mid L(h) = 0\}.$  We see from Corollary 4.3
of \cite{hufo} that for any $c_0, c_1, \ldots, c_{n-1} \in k$, there
exists a unique $y \in V$ such that $\pi_j(y) = c_j$ for $j = 0,
1, \ldots, n-1$.

\bprop \label{dim-sol}
    Let $h_0, h_1, \ldots, h_{n-1} \in Hk$, and let $L$ be the
linear homogeneous differential operator on $Hk$ defined for any $h
\in Hk$ by
$$L(h) =\partial^{n}(h)+ \sum_{i=0}^{n - 1} h_{i}\partial^{i}(h)
.$$
Then $V$ is an $n$-dimensional $k$-vector space.
\eprop

\bpf
    Since $L: Hk \to Hk$ is a $k$-linear operator, it
is clear that $ V = \ker(L)$ is a $k$-vector space, so it remains
to prove that $\dim_{k} V = n$. To see this, we define a mapping
$T: k^{n} \to V$ as follows: If $\underline{a} = (a_{1}, \ldots ,
a_{n}) \in k^{n},$ then $T(\underline{a})$ is the unique solution in
$Hk$ to $L(h) = 0$ such that $\pi_{i}(T(\underline{a})) = a_{i+1}$
for $i = 0, \ldots , n-1$ by~\cite[Corollary 4.3]{hufo}. It is clear
that $T$ is a $k$-vector space isomorphism, from which the result
follows. \epf

It follows that $Hk$ has the following ``completeness'' property:
Any $n^{th}$ order monic linear homogeneous ordinary differential
equation with coefficients in $Hk$ has a complete set of $n$
linearly independent solutions in $Hk$.

This can be done more generally as follows. Let $A$ denote any
commutative ring with identity, let $h_0, h_1, \ldots, h_{n-1} \in
HA$ and $c_0, c_1, \ldots, c_{n-1} \in A$.  As before, consider the
linear homogeneous differential operator $L$ defined on $HA$ for any
$h \in HA$ by $L(h) = \partial^{n}(h)+\sum_{i=0}^{n-1} h_{i}\partial^{i}(h)
$.  We know from Corollary 4.3 of \cite{hufo} that
for any $c_0, c_1, \ldots, c_{n-1} \in A$, there is a unique
solution $y \in HA$ to $L(h) = 0$ such that $\pi_i(y) = c_i$ for
each $i = 0, 1, \ldots, n-1$. We now give a constructive method for
finding solutions to $L(h) = 0$ in $HA$.

\bprop \label{recurrence}
Let $A$ be a commutative ring with identity, let $h_i \in HA$ and
let $c_i \in A$ for $i = 0 , \ldots, n-1$. Let $L$ be the linear
homogeneous differential operator defined on $HA$ for any $h \in HA$
by $$L(h) = \partial^{n}(h)+\sum_{i=0}^{n-1} h_{i}\partial^{i}(h)
.$$  The unique solution $y \in HA$ to $L(h) = 0$
such that $\pi_i(y) = c_i$ for each $i = 0, 1, \ldots, n-1$ is
given by
\begin{equation}
\pi_i(y) = c_i, i = 0, 1, \ldots, n-1
\end{equation}
and
\begin{equation}
\pi_{n+m}(y) = -\sum_{i=0}^{n-1}\sum_{j=0}^{m}{m \choose j}
\pi_j(h_i)\pi_{m-j+i}(y), m \in \mathbf{N}.
\end{equation}
\eprop

\bpf
    Clearly $y \in HA$ given by the above prescription is unique,
and $y$ satisfies the initial conditions $\pi_i(y) = c_i, i = 0,
1, \ldots, n-1$ by definition, so we must only show that $L(y) =
0$.  This means we must show that for each $r \in \mathbf{N}$,
$\pi_{r}(L(y)) = 0$.  Now we have
\begin{eqnarray*}
\pi_r(L(y))
  & = &
\sum_{i=0}^{n-1} \pi_r(h_{i}\partial^{i}(y)) +
\pi_r(\partial^{n}(y)) \\
  & = &
\sum_{i=0}^{n-1}\sum_{j=0}^{r}{r \choose j}
\pi_j(h_i)\pi_{r-j}(\partial^i(y)) +
\pi_r(\partial^{n}(y)) \\
  & = &
\sum_{i=0}^{n-1}\sum_{j=0}^{r}{r \choose j}
\pi_j(h_i)\pi_{r-j+i}(y) + \pi_{r+n}(y) \\
  & = &
\sum_{i=0}^{n-1}\sum_{j=0}^{r}{r \choose j}
\pi_j(h_i)\pi_{r-j+i}(y) + (-\sum_{i=0}^{n-1}\sum_{j=0}^{r}{r
\choose j} \pi_j(h_i)\pi_{r-j+i}(y)) \\
  & = &
0.\end{eqnarray*}\epf

The following corollary gives a very simple description of the
solutions in the case that the coefficients of the equation are
constants.

\bco \label{linear recursive}
Let $A$ be a commutative ring with identity, let $a_0, \ldots,
a_{n-1} \in A$ and let $c_0, \ldots, c_{n-1} \in A$. Let $L$ be the
linear homogeneous differential operator defined on $HA$ for any $h
\in HA$ by $$L(h) = \partial^{n}(h)+\sum_{i=0}^{n-1} a_{i}\partial^{i}(h)
.$$  The unique solution $y \in HA$ to $L(h) = 0$
such that $\pi_i(y) = c_i$ for each $i = 0, 1, \ldots, n-1$ is given
by
\begin{equation}
\pi_i(y) = c_i, i = 0, 1, \ldots, n-1
\end{equation}
and
\begin{equation}
\pi_{n+m}(y) = -\sum_{i=0}^{n-1}a_i\pi_{m+i}(y), m \in \mathbf{N}.
\end{equation}
or more simply,
\begin{equation}
y_{n+m} = -\sum_{i=0}^{n-1}a_iy_{m+i}, m \in \mathbf{N}.
\end{equation}
\eco

\bpf Since the $h_i = a_i$ are constants, we have
$\pi_j(h_i) = a_i$ if $j = 0$  and $\pi_j(h_i) = 0$ if $j \geq 1$.
Therefore the only nonzero term in the inner sum is the $j = 0$
term.  From this the result follows.\epf

Corollary \ref{linear recursive} shows that, in the case of constant
coefficients, the solutions $y\in Hk$ to $L(Y) = 0$ are
\textsl{linearly recursive sequences}.

\section{Linear Homogeneous Differential Equations with Constant Coefficients}

As before, let $k$ be a field of any characteristic  and let $Hk$
be the differential ring of Hurwitz series over $k$.  For any
$\beta \in k$, the element $\exp(\beta) = (1, \beta, \beta^{2},
\ldots, \beta^{n}, \ldots) \in Hk$ is called the
\textsl{exponential} of $\beta$.  Note that for any $c \in k$,
$c\exp(\beta)$ is the unique solution in $Hk$ to the differential
equation $\partial(y) - \beta y =0$ with initial condition $y(0) =
c$.  The following result is immediate.

\bl \label{exp}
    Let $\alpha, \beta \in k$.  Then
\begin{enumerate}
  \item $\exp(\alpha + \beta) = \exp(\alpha) \exp(\beta)$
  \item $\exp(0) = 1$
  \item For each $\beta \in k$, $\exp(\beta)$ is invertible in
  $Hk$, and $\exp(-\beta) = \exp(\beta)^{-1}$.
\end{enumerate}
\el

From~\cite{hufo} we recall the divided powers $x^{[i]}$ in $Hk$, for
$i \in \mathbf{N}$, defined by $x^{[i]} = (\delta_n^i)$, so that
$$x^{[0]} = 1_{Hk}, \hspace{.1in} x^{[1]} = x = (0, 1, 0, 0, \ldots,
0, \ldots), \hspace{.1in} x^{[2]} = (0, 0, 1, 0, \ldots),$$ etc.
Using the natural topology on $Hk$ and the divided powers $x^{[i]}$,
we have $\exp(\beta) = \sum_{i=0}^{\infty} \beta^{i} x^{[i]} =
\sum_{i=0}^{\infty} (\beta x)^{[i]}$.  We will denote $\exp(\beta)$ by $e^{\beta x}$.

For simplicity, assume that $k$ is algebraically closed, although we
shall see later that this assumption won't be necessary. Let $V$ be a $k$-subspace of $Hk$ that is closed under the  derivation $\partial$
which we denote by $'$. We denote the group of all k-differential automorphisms of $V$ by $G(V|k)$. That is, $G(V|k):=\{\sigma \in
\textrm{Aut}_k V \ | \  \sigma(v)' =
\sigma(v')\ \textrm{for all}\  v \in V\}$. We will sometimes denote $G(V|k)$ by $G(L|k)$ if $V$ is the full set of solutions
of a linear homogeneous differential equation $L(Y)=0$.

\subsection{Computing the group $G(V|k)$.}\label{group} Let $y \in Hk$ be a non-zero solution of the equation $Y' = \alpha Y$ where $\alpha\in k$ and let $\pi_0(y)=c$.
Then from Corollary 4.3 of \cite{hufo}, it follows that $y=ce^{\alpha x}$.  Thus $V = \{ce^{\alpha x}| c \in
k\}$ forms a full set of solutions of the equation $Y' = \alpha Y$. Let $\sigma\in G(V|k)$ and note that
 $\sigma (e^{\alpha x})'=\alpha \sigma (e^{\alpha x})$. Thus
$\sigma(e^{\alpha x})$ is also a solution of $Y' = \alpha Y$. Therefore
there is some $c_{\sigma} \in k^*$ such that $\sigma(e^{\alpha x}) =
c_{\sigma}e^{\alpha x}$.  Also note that $\sigma\tau(e^{\alpha x}) =
c_{\tau}\sigma(e^{\alpha x}) = c_{\tau}c_{\sigma}e^{\alpha x} = c_{\sigma}c_{\tau}e^{\alpha x} =
\tau\sigma(e^{\alpha x})$.  Thus $G(V|k) \hookrightarrow  (k^*, \times)$ by
$\sigma \mapsto c_{\sigma}$ is an injective group homomorphism.
Moreover, for any $c \in k^*$, we may define $\sigma_c: V
\rightarrow V$ by $\sigma_c(e^{\alpha x}) = ce^{\alpha x}$.  It is clear that $\sigma_c$
is a $k$-differential automorphism of $V$. Thus $G(V|k) \cong (k^*, \times)$.

More generally, for $\alpha \in k,$ consider the $k$-subspace $V_t$ of $Hk$ defined by $$V_t =
\textrm{Span}_k\{z_0,z_1,\cdots,z_t\},$$  where $z_j = x^{[j]}e^{\alpha x}$ for $j = 0,1, \cdots, t$. It can be shown that $\{z_0,z_1,\cdots,z_t\}$ is linearly independent over $k$, since w($z_0,z_1,\cdots,z_t$)$= z^{t+1}_0=e^{(t+1)\al x}$, which is invertible in $Hk$. One can construct a linear
homogeneous differential equation over $k$ whose full set of solutions equals $V_t$, namely
\begin{equation}\label{wroneqn}L_t(Y) = \frac{w(Y,z_0,z_1,\cdots,z_t)}{w(z_0,z_1,\cdots,z_t)}.\end{equation}

\bt\label{repeatedroots} Let $\al\in k$ and let $z_j = x^{[j]}e^{\alpha x}$ for $j = 0,1, \cdots, t$.  Let $V_t =
\textrm{Span}_k\{z_0,z_1,\cdots,z_t\}.$
Let $Z_t:=\begin{pmatrix}
z_0, &z_1, & \cdots, &z_t
\end{pmatrix}$ and let $I_{t}$ be the identity matrix of order $t$. Let $\ur_t:=\begin{pmatrix}0 & I_{t}\\0 & 0\end{pmatrix}\in M(t+1,k)$ if $t\geq 1$, and $\ur_0=0\in k$.
Then
\begin{enumerate}

\item $Z_t'=Z_t(\alpha I_{t+1}+\ur_t)$.\\

\item
$C_k(I_{t+1}+\ur_t)=C_k(\ur_t)=Span_{k}\{I_t, \ur_t,\ur^2_t,\cdots,\ur^{t-1}_t\}.$\\

\item  Under the basis $z_0,z_1,\cdots, z_t$, the group

\begin{equation} G(V_t|k)\cong
\left\{
  \begin{array}{ll}
    C_k(\ur_t), & \hbox{if $\al\neq 0$;} \\
    C_k(\ur_t)\cap U(t+1,k), & \hbox{if $\al=0$;}
      \end{array}
\right.\end{equation}
where $U(t+1,k)$ is the group of all upper triangular matrices in $GL(t+1,k)$ with $1$ on the main diagonal.
\end{enumerate}

\et

\bpf

Since $z'_i=\alpha z_i+z_{i-1}$ for all $i\geq 1$ and $z'_0=\al z_0$, it follows that $Z_t'=Z_t(\alpha I_{t+1}+\ur_t)$.

A straightforward computation proves (b).

Let $\s\in G(V_t|k)$. Since  $\s(z_j)\in V_t$, there is an element $C_\s\in G$ such that $\s(Z_t)=Z_tC_{\s}$. Therefore
\begin{align*}\s(Z'_t)&=\s(Z_t)(\alpha I_{t+1}+\ur_t)\\
&=Z_tC_{\s}(\alpha I_{t+1}+\ur_t)\end{align*}
On the other hand, $\s(Z_t)'=(Z_tC_{\s})'=$ $Z_t(\alpha I_{t+1}+\ur_t)C_{\s}$. Since $\s(Z_t)'=\s(Z'_t)$, we obtain that $(\alpha I_{t+1}+\ur_t)C_{\s}$
$=C_{\s}(\alpha I_{t+1}+\ur_t)$, which is true if and only if  $C_{\s}\ur_t=\ur_t C_{\s}$.  Thus there is an injective group homomorphism $\phi: G(V_t|k)\to C_k(\ur_t)$ given by
$\phi(\s)= C_\s$. Moreover, if $\al=0$ then $z_0=1$ and therefore $\s(z_0)=z_0$. It then follows that $\phi(G(V_t|k))\subseteq C_k(\ur_t)\cap U(t+1,k)$ if $\al= 0$.

To prove that $\phi$ is surjective in each of the cases $\al=0$ and $\al\neq 0$, we first note that for any constants $c_0,c_1,\cdots,c_t\in k$ and $0\leq j\leq t$
\begin{align}\notag(\sum^j_{i=0}c_{j-i}z_i)'&=c_jz'_0+\sum^j_{i=1}c_{j-i}z'_i\\
\notag &=\al c_jz_0+\sum^j_{i=1}c_{j-i}(\al z_i+z_{i-1}) \\
\label{diffauto} &=\sum^j_{i=0}\al c_{j-i}z_i+\sum^j_{i=1}c_{j-i}z_{i-1}.
\end{align}

Let $C\in C_k(\ur_t)$  when $\al\neq 0$ and $C\in C_k(\ur_t)\cap U(t+1,k)$ when $\al=0$. Then we may define an automorphism of the $k$-vector space $V_t$ such that $\s_C(Z)=ZC$, that is $\s_C(z_j)=\sum^j_{i=0}c_{j-i}z_i$. Now from Equation \ref{diffauto} it can be seen that $\s_C$ is a $k$-differential automorphism of $V$. Thus $G(V_t|k)\cong C_k(\ur_t)$ or $G(V_t|k)\cong C_k(\ur_t)\cap U(t+1,k)$ depending on whether $\al\neq 0$ or $\al=0$ respectively.\epf

The following theorem is an immediate consequence of Theorem \ref{repeatedroots}
\bt\label{groupforconstcoeff}
Let $\al_1,\al_2,\cdots,\al_r\in k$ be distinct elements. For $t=1,2,\cdots,r$, let $V_t =
\textrm{Span}_k\{z_{j,t}|\ 0\leq j\leq m_t\}$ and $V=\bigoplus^r_{t=1}V_t$,  where $z_{j,t} = x^{[j]}e^{\al_t x}$.
Let $Z:=\begin{pmatrix}
z_{0,1},& \cdots, &z_{m_1,1},&\cdots,& z_{0,r},&\cdots,&z_{m_r,r}
\end{pmatrix}$, $I_{m_t}$ be the identity matrix of order $m_t$ and $\ur_t:=\begin{pmatrix}0 & I_{m_t}\\0 & 0\end{pmatrix}$ $\in M(m_t+1,k)$.
Then

\begin{enumerate}

\item $Z'=Z(S+\ur)$,  where $$S=\begin{pmatrix}
\al_1I_{m_1+1} & 0 & \cdots&0 \\
0 & \al_2I_{m_2+1} &  \cdots& 0\\
\vdots & \vdots& \ddots & \vdots\\
0 & 0 &  0 & \al_rI_{m_r+1}
\end{pmatrix},\quad \quad  \ur=\begin{pmatrix}
\ur_1 & 0 & \cdots&0 \\
0 & \ur_2 &  \cdots& 0\\
\vdots & \vdots& \ddots & \vdots\\
0 & 0 &  0 & \ur_r
\end{pmatrix}.$$\\

\item We have the following group isomorphism
\begin{equation}\label{gpiso}G(V|k)\cong \bigoplus^r_{t=1}G(V_t|k). \end{equation}

\end{enumerate}

\et
Since $z_{0,1},\cdots, z_{m_1,1},\cdots,z_{0,r},\cdots, z_{m_r,r}$ is a basis for $V$, each of the groups $G(V_t|k)$ can be computed using Theorem \ref{repeatedroots}.

\subsection{Non-algebraically closed fields}

Let $k$ be a field  and let $\bar{k}$ denote its algebraic closure. Let $a_0,\ldots,a_{n-1} \in k$ and consider the monic linear homogeneous differential operator $L(y) = y^{(n)}+\sum_{i=0}^{n-1} a_{i}y^{(i)}.$  Then from Proposition \ref{dim-sol}, we know that there are $k$-linearly independent elements $y_1,y_2,\cdots,y_n\in Hk$ such that $L(y_i)=0$ for each $i$. Let $Y:=(y_1,\cdots, y_n)$ and note  that $y'_i$ is also a solution of $L(y)=0$ for each $i$. Thus, there is a matrix $B\in M(n,k)$ such that \begin{equation}\label{matrix eqn}Y'=YB.\end{equation}

Considering the differential operator $L(y)$ over the field $\bar{k}$, one can factor the characteristic polynomial of the operator $L(y)$. Let $\al_1,\al_2,\cdots,\al_r$ be the distinct roots of the characteristic polynomial in $\bar{k}$. Let $Z, S, \ur$ be as in Theorem \ref{groupforconstcoeff}.
  Then it can be shown that  $V(Z,\bar{k}) :=$ Span$_{\bar{k}}\{z_1,\cdots,z_n\}$ is the set of all solutions of $L(y)=0$ in $H\bar{k}$. Let  $V(Y,k):=$ Span$_k\{y_1,\cdots,y_n\}$ and note that  $V(Y,k)\subset V(Z,\bar{k})$ and since $w(y_1,\cdots,y_n)\neq 0$, $y_1,y_2,\cdots,y_n$ remain linearly independent over $\bar{k}$. Thus  $V(Z, \bar{k})=V(Y,\bar{k})$. Let $\phi:V(Y, \bar{k})\to V(Z,\bar{k})$ be a map of $\bar{k}$-vector spaces such that $\phi$ maps the ordered basis $Y$ to the ordered basis $Z$. Then there is a matrix $T_{\phi}\in GL(n,\bar{k})$ such that $YT_{\phi}=Z$. Applying $\phi$ to the Equation \ref{matrix eqn}, we obtain
 \begin{align}
\notag Z'&= Y'T_\phi\\
\notag &=YBT_\phi\\
&=ZT^{-1}_{\phi}BT_\phi.
  \end{align}
Thus we see that $T^{-1}_\phi BT_\phi=S+\ur$. In particular, $S+\ur$ is the Jordan normal form of $B$. From the above discussion, we derive the following theorem.

\bprop \label{compcent} Let $k$ be a field and let $L(y) = y^{(n)}+\sum_{i=0}^{n-1} a_{i}y^{(i)},$ where $a_i\in k$ for each $i$. Let $\al_1,\cdots,\al_r$ be the distinct roots of the characteristic polynomial of $L(y)=0$ in $\bar{k}$. Let $Z, S,\ur$ be as in Theorem \ref{groupforconstcoeff} and let $A:=S+\ur$. Let $y_1,y_2,\cdots,y_n\in Hk$ be $k$-linearly independent solutions of $L(y)=0$ in $Hk$, $Y:=(y_1,\cdots, y_n)$ and let $Y'=YB$ for some $B\in M(n,k)$.
. Then \begin{enumerate}
\item Span$_k\{B^i|0\leq i\leq n-1\}=$ Span$_k\{B^i|0\leq i\leq \infty\}$, where $B^0:=I$.

\item if $a_0\neq 0$ then $$G(V|k)\cong C_{k}(B)=\text{Span}_{k}\{B^i|0\leq i\leq n-1\}.$$

\end{enumerate}
\eprop

\bpf
Since $Y'=YB$, we have $B^n=-\sum_{i=0}^{n-1} a_{i}B^i.$ To prove (1), it is enough to show that $\{I, B, \cdots,B^{n-1}\}$ is linearly independent over $\bar{k}$. Suppose that $b_0,\cdots,b_{n-1}\in \bar{k}$ and $\sum_{i=0}^{n-1} b_{i}B^i=0$. Then we have $\sum_{i=0}^{n-1} b_{i}YB^i=0$, which implies $\sum_{i=0}^{n-1} b_{i}Y^{(i)}=0$. Let $G(y):=\sum_{i=0}^{n-1} b_{i}y^{(i)}$ and note that $G(y_j)=0$ for each $j=1,2,\cdots,n.$ Since the order of $G(y)$ is less than $n$, we obtain that $b_i=0$ for all $i$. Thus (1) is proved.

It suffices to consider the case when $r=1$. We know from Theorem \ref{repeatedroots} that $C_{\bar{k}}(A)=$Span$_{\bar{k}}\{I, \ur,\ur^2,\cdots,\ur^{n-1}\}$. Let $T\in GL(n,\bar{k})$ such that $YT=Z$. Then since $C_{\bar{k}}(B)=TC_{\bar{k}}(A)T^{-1}$, we see that $C_{\bar{k}}(B)$ is a $\bar{k}$-vector space of dimension $n$. Now $B\in C_{\bar{k}}(B)$ will imply that $C_{\bar{k}}(B)=$ Span$_{\bar{k}}\{I, B, \cdots,$ $B^{n-1}\}$. Since $\{I, B, \cdots,$ $B^{n-1}\}$ is $\bar{k}-$linearly independent, it follows that $C_{k}(B)=C_{\bar{k}}(B)\cap GL(n,k)=$ Span$_{k}\{I, B, \cdots,$ $B^{n-1}\}$.   \epf

\bt \label{nonalgcase}
Let $k$ be a field, $\bar{k}$ be its algebraic closure and let $L(y) = y^{(n)}+\sum_{i=0}^{n-1} a_{i}y^{(i)},$ where $a_i\in k$ for each $i$. Let $\al_1,\cdots,\al_r$ be the distinct roots of the characteristic polynomial of $L(y)=0$ in $\bar{k}$. Let $Z, S,\ur$ be as in Theorem \ref{groupforconstcoeff} and let $A:=S+\ur$. Let $y_1,y_2,\cdots,y_n\in Hk$ be $k$-linearly independent solutions of $L(y)=0$ in $Hk$, $Y:=(y_1,\cdots, y_n)$ and let $Y'=YB$ for some $B\in M(n,k)$.

Then, \begin{enumerate}
\item if $a_0\neq 0$, \begin{equation}\label{group formula2} G(V|k)\cong C_{k}(B)=\text{Span}_{k}\{B^i|0\leq i\leq n-1\}, \end{equation}

\item if $a_0=0$, \begin{equation}\label{group formula1}G(V|k)\cong GL(n,k) \cap T\begin{pmatrix}
\g_1 & 0 & \cdots&0 \\
0 & \g_2 &  \cdots& 0\\
\vdots & \vdots& \ddots & \vdots\\
0 & 0 &  0 & \g_r
\end{pmatrix}T^{-1},\end{equation}
where $T\in GL(n,\bar{k})$ such that $YT=Z$, $\g_i= C_{\bar{k}}(\ur_i)\cap U(m_i+1,\bar{k})$ for at most one $i$ and in that case $\al_i=0$ and $\g_t=C_{\bar{k}}(\ur_t)$ for all other $t$.
\end{enumerate}
\et
\bpf
Note that  $C_{\bar{k}}(B)=T C_{\bar{k}}(A) T^{-1}$ and therefore $C_{k}(B)=GL(n,k)\cap TC_{\bar{k}}(A) T^{-1}$. Now the rest of the proof follows from Theorem \ref{groupforconstcoeff}. \epf
{\bf Remark.} From Theorem \ref{repeatedroots} (2), we see that $C_{\bar{k}}(u_t)$ is a commutative linear algebraic group for each $t, 1\leq t\leq r$. Then it follows that $G(V|k)$ is a commutative linear algebraic group as well. Also, from Theorem \ref{nonalgcase}, we see that the condition that $k$ be algebraically closed is not needed.


\section{Examples}

Let $k$ be a field of any characteristic. In the following examples, we will compute the group $G(V|k)$.

{\bf Example 1.}
Consider the second order operator $L(y)=y''$.  Let $Y=(1,x)$ and note that $V= {\rm Span}_kY$ consists of
all solutions of the equation $y''=0$. Also note that $$Y'=YA,$$ where $A=\begin{pmatrix}0 & 1\\0 & 0\end{pmatrix}.$
Applying Theorem \ref{groupforconstcoeff}, we obtain that the group $G(V|k)\cong U(2,k).$

{\bf Example 2.}
Consider the differential equation $$L(Y)=Y''-Y'-Y.$$ Let $Y =(y_1,y_2)$,  where $y_1=(1,0,1,1,2,3,\cdots)$ and  $y_2=(0,1,1,2,3,5,\cdots)$ are Fibonacci sequences. Then it can be seen that $V= {\rm Span}_k\{y_1,y_2\}$ consists of all solutions of $L(y)=0$. Since $a_0=-1$, from Equation \ref{group formula2}, we obtain that $$G(V|k)\cong C_{k}(B)$$ with respect to the basis $Y$, where $B=\begin{pmatrix}0 & 1\\ 1 & 1\end{pmatrix}$.
One can directly compute the centralizer of $B$ and obtain \begin{align*}C_{k}(B)&=\left\{\begin{pmatrix}\al & \beta\\\beta & \al+\beta\end{pmatrix}\in GL(2,k)|\al,\beta\in k \right\}\\
&= \textrm{Span}_k\{I, B\}\cap GL(2,k).\end{align*}

{\bf A note on initial conditions.} Let $a_0,\ldots,a_{n-1} \in k$ and consider the monic linear homogeneous differential operator $L(y) = y^{(n)}+\sum_{i=0}^{n-1} a_{i}y^{(i)}, y\in Hk.$ Let $y_1,y_2,\cdots,y_n\in Hk$ be  $k$-linearly independent elements such that $L(y_i)=0$ and that $\pi_{i-1}(y_j)=\delta^j_i$ for each $i,j=1,2,\cdots,n$. Let $Z:=(z_1,z_2,\cdots,z_n)$, where $\{z_1,\cdots,z_n\}$ is set of linearly independent solutions of $L(y)=0$ in $H\bar{k}$. Let $Y=(y_1,y_2,\cdots,y_n)$ and $Z=(z_1,z_2,\cdots,z_n)$. Then it follows, from the uniqueness of solutions subject to initial conditions, that for \begin{equation}\label{choiceofT}T:=\begin{pmatrix}
\pi_0(z_1) & \pi_0(z_2) & \cdots&\pi_0(z_n) \\
\pi_1(z_1) & \pi_1(z_2) & \cdots&\pi_1(z_n)\\
\vdots & \vdots& \vdots & \vdots\\
\pi_{n-1}(z_1) & \pi_{n-1}(z_2) & \cdots&\pi_{n-1}(z_n)
\end{pmatrix}\end{equation}
we have $YT=Z$. This observation along with Equation \ref{group formula1}, enables us to compute the group  $G(V|k)$ with respect to the basis $Y$.

{\bf Example 3.} Consider the operator $L(Y)=Y'''-3Y''+3Y'-Y$ and let $Y=(y_1,y_2,y_3)$ be linearly independent solutions of $L(Y)=0$ with initial conditions $\pi_{i-1}(y_j)=\delta^j_i$ for each $i,j=1,2$ and $3$. Let $B=\begin{pmatrix}
0 & 1 & 0\\
0 & 0 & 1\\
1 & -3 & 3
\end{pmatrix}$ and  note that $Y'=YB$. From Proposition \ref{compcent}, it follows that
\begin{align}
\notag G(V|k)&\cong C_{k}(B)\\
\label{span eqn}&=\text{Span}_{\bar{k}}\{I,B,B^2\}\cap GL(3,k)
\end{align}

It is  also possible to realize the group as a full set solutions for a system of linear equations over $k$. Note that $$ C_{\bar{k}}(\ur_2)= \left\{\begin{pmatrix}
a & b & c\\
0 & a & b\\
0 & 0 & a
\end{pmatrix} \in GL(3,\bar{k})\right\}$$ and that $YT=Z$ for $T=\begin{pmatrix}
1 & 0 & 0\\
1 & 1 & 0\\
1 & 2 & 1
\end{pmatrix}$. Thus we have

\begin{align*}G(V|k)&\cong C_{k}(B)= TC_{\bar{k}}(\ur_2)T^{-1}\cap GL(3,k)\\
&= \left\{\begin{pmatrix}
a-b+c & b-2c &c\\
c & -b+a-2c & c+b\\
c+b & -3b-2c & c+2b+a
\end{pmatrix} | a,b,c\in k, a\neq 0\right\}.\end{align*}
If char($k$)=2 we note that
$$G(V|k)= \left\{\begin{pmatrix}
a+b+c & b &c\\
c & a+b & b+c\\
b+c & b& a+c
\end{pmatrix} | a,b,c\in k, a\neq 0\right\}.$$


\begin{thebibliography}{99}



\bibitem{onthri} W. F. Keigher, On the ring of Hurwitz series,
{\em Comm. Algebra} {\bf 25} (1997), 1845-1859.

\bibitem{hufo} W. F. Keigher and F. L. Pritchard, Hurwitz series as formal
functions, {\em J. Pure Appl. Algebra} {\bf 146} (2000), 291-304.



\end{thebibliography}
\end{document}